\documentclass[12pt]{article}

\if@twoside
\else
\fi \marginparwidth 0pt \marginparsep 0pt \topmargin -1.4cm
\usepackage{amsmath,amsthm,amsfonts,amssymb}
\usepackage[latin1]{inputenc}
\usepackage{makeidx}
\usepackage[french]{babel}
\usepackage{newlfont}
\usepackage{showidx}
\usepackage[dvips]{graphicx}

\newtheorem{theo}{Théorème}
\newtheorem{prop}{Proposition}

\newtheorem{lem}{Lemme}

\def\g{\mathfrak{g}}
\def\G{\Gamma}
\def\ov{\overline}
\def\tr{\mathrm{tr}}
\def\ad{\mathrm{ad}}
\def\p{\mathfrak{p}}
\def\k{\mathfrak{k}}
\def\Exp{\mathrm{Exp}}
\newcommand{\R}{\mathrm{ I\! R}}
\newcommand{\C}{\mathrm{ I\!\!\!C}}

\begin{document}

\bibliographystyle{alpha}
\author{Charles Torossian\footnote{UMR 8553 du CNRS, DMA-ENS, Ecole Normale Supérieure, 45 rue d'Ulm 75230 Paris
cedex 05, Charles.Torossian@ens.fr}}
\title{ Méthodes de Kashiwara-Vergne-Rouvière pour certains espaces symétriques}
\date{}
\maketitle

\begin{abstract} Cet article est une suite de notre article
\cite{To}. En utilisant une déformation  "à la Kontsevich" de la
formule de Campbell-Hausdorff pour les espaces symétriques, on
retrouve les résultats de Rouvière \cite{rou86}, sur la convolution
des distributions invariantes, dans le cas des espaces symétriques
résolubles et "très symétriques".
\\

This paper follows our previous work \cite{To}. We study the case of
symmetric spaces. We recover, by using a Kontsevich's deformation of
the Baker-Campbell-Hausdorff formula, Rouvière's results
\cite{rou86}, on the convolution of invariant distributions, for
solvable symmetric spaces and "very symmetric spaces".
\end{abstract}
\section*{Introduction}
Dans notre article précédent \cite{To}, suivant les idées élaborées
par Kontsevich \cite{kont} pour la quantification formelle des
variétés de Poisson et déjà utilisées dans \cite{ADS} et \cite{AST},
nous avons montré comment résoudre certaines conjectures sur la
formule de Campbell-Hausdorff (BCH in English). Plus précisément,
nous avons construit une défor\-mation de la formule BCH, qui
vérifie une équation différentielle analogue à celle que l'on trouve
dans l'article  de Kashiwara-Vergne \cite{KV}. Ces équations ont la
même utilité que celles de l'article de Kashiwara-Vergne. En effet
elles permettent de montrer, comme dans \cite{ADS} et \cite{AST},
que l'application exponentielle modifiée par la racine carrée du
jacobien, transporte la convolution des distributions invariantes
sur les algèbres de Lie (\cite{Du},\cite{KV}, \cite{Ve},\cite{Mo}).
\\

Dans cet article nous considérons le cas des espaces symétriques.
\\

Dans notre mémoire d'Habilitation \cite{hab} nous avons par ailleurs
abordé certains problèmes d'analyse harmonique pour les espaces
symétriques. Nous avons notamment expliqué en quoi le comportement
des espaces symétriques est à la fois proche et éloigné du cas des
algèbres de Lie. On peut dire en résumant grossièrement, qu'il
existe une catégorie d'espaces symétriques pour laquelle les
méthodes des groupes fonctionnent bien (cas nilpotent, résoluble,
$"G_{\C}/G_{\R}"$, cas "très symétrique"). Malheureusement le cas
général des espaces symétriques ne peut pas être traité par des
méthodes s'inspirant simplement du cas des groupes (voir la fin de
cet article pour  quelques remarques).
\\

Dans cet article nous allons déformer la formule de
Campbell-Hausdorff pour les espaces symétriques généraux.
\\

Nous montrerons que cette déformation vérifie, dans les "bons" cas
évoqués plus haut (résoluble ou très-symétrique), des équations
analogues au cas des groupes. On en déduira que l'application
exponentielle modifiée par la racine carrée du déterminant jacobien,
transporte la convolution des germes de distributions invariantes.
Nous retrouvons ainsi les résultats de Rouvière \cite{rou86} du cas
résoluble.  Notre démonstration fournit une preuve dans le cas très
symétrique ce qui, à notre connaissance, est un résultat nouveau.
\\

\noindent {\bf Remerciements:} L'auteur remercie M. Pevzner et le
département de mathéma\-tiques l'Université Libre de Bruxelles (ULB)
pour leur invitation (fin janvier 2002). L'auteur y a finalis\'e les
résultats de cet article. Nous remercions aussi les organisateurs du
colloque Carmona pour leur invitation. Nous avons eu le plaisir d'y
pr\'esenter les r\'esultats de cet article.

\section{Rappels}

\subsection{Rappels sur la formule de Kontsevich}

La formule de Kontsevich associe à toute structure de Poisson
régulière sur $\R^d$ un étoile-produit formel associatif. Ce n'est
qu'un cas particulier du théorème de formalité d\'emontr\'e par
Kontsevich dans \cite{kont} théorème $6.4$.

Lorsque   $f$ et $g$ sont deux fonctions régulières sur $\R^{d}$ et
$\alpha$ un deux-vecteur de Poisson régulier sur $\R^d$, Kontsevich
écrit dans \cite{kont} paragraphe $2$, la formule suivante :

\begin{equation}\label{formuleK}
f\star g=fg +\sum _{n=1}^{\infty}\frac{h^{n}}{n!}\sum_{\substack{\G
\in G_{n,2}\\ \Gamma {\mathrm admissible}}}w_{\Gamma} B_{\G}(f,g).
\end{equation}

Cette formule munit l'espace des fonctions régulières sur $\R^d$
d'une structure associative formelle. Dans cette formule $h$ est un
paramètre formel, $G_{n,2}$ désigne l'ensemble des graphes étiquetés
avec $n$ points de première espèce et $2$ points de seconde espèce,
$\G$ est un graphe dit admissible parmi les graphes de $G_{n,2}$,
$w_{\Gamma}$ est un c\oe fficient obtenu par intégration sur un
espace de configurations d'une forme différentielle dépendant de
$\G$ et $B_{\G}$ est un opérateur bidifférentiel construit à partir
de $\G$.

Par souci de clarté et d'autonomie du texte nous allons préciser
maintenant chacun des termes de cette formule et appliquer dans un
second temps cette formule dans le cas des structures de Poisson
linéaires c'est-à-dire dans le cas du dual des algèbres de Lie.

\subsubsection{Espaces de configurations}

On note par $\mathrm{Conf}_{n,m}$ l'espace des configurations de $n$
points distincts dans le demi-plan de Poincaré (ce sont les points
aériens) et de $m$ points  distincts sur la droite réelle (ce sont
les points terrestres). Le groupe:
$$G^{(1)} =\{z\mapsto az+b \hbox{ avec }(a,b)\in\R^*_+\times\R\}$$
agit librement sur $\mathrm{Conf}_{n,m}$. Le quotient~:
$$C_{n,m}=\mathrm{Conf}_{n,m}/G^{(1)}$$
est une vari\'et\'e de dimension $2n+m-2$. Compte tenu de l'action
de ce groupe sur les points terrestres, on peut identifier deux des
points terrestres aux points $0$ et $1$ (à condition que l'on ait
$m\geq 2$, sinon on peut identifier un des points aériens au
complexe $i$). Dans \cite{kont} paragraphe $5.1$, Kontsevich
construit des compactifications de ces variétés notées
$\overline{C}_{n,m}$. Ce sont des variétés à coins de dimension
$2n+m-2$. Ces variétés ne sont pas connexes pour $m\geq 2$. On
notera par $\overline{C}^{+}_{n,m}$ la composante connexe qui
contient les configurations où les points terrestres sont ordonnés
dans l'ordre croissant (i.e.  on  a $\ov{1}< \ov{2}<\cdots<\ov{m}$).
On introduit de manière analogue des variétés de configurations de
$n$ points dans le plan complexe modulo l'action du groupe
$G^{(2)}=\{z\mapsto az+b \hbox{ avec }(a,b)\in\R_+^*\times\C \}$. On
les note $C_n$, ce sont des variétés de dimension $2n-3$. On note
$\overline{C_n}$ les compactifications associées (\cite{kont}
paragraphe $5.1$).

Ces variétés sont stratifiées et chaque strate est décrite par un
arbre. En termes géométriques, les strates sont obtenues par
concentrations it\'er\'ees de points en des amas (\cite{kont}
paragraphe 5.2).

\subsubsection{Graphes admissibles}\label{defgrapheadmissible}
La notion de graphes  admissibles  est maintenant bien établie dans
la littéra\-ture. On d\'esigne par $G_{n,m}$ l'ensemble des graphes
\'etiquet\'es et orient\'es (les arêtes sont orientées) ayant $n$
sommets du première espèce (sommets aériens) et $m$ sommets du
deuxième espèce (sommets terrestres). Par graphe \'etiquet\'e on
entend un graphe $\Gamma$ muni d'un ordre total sur l'ensemble
$E_\Gamma$ de ses ar\^etes, compatible avec l'ordre des sommets.

Les graphes qui vont intervenir dans la formule de Kontsevich sont
dans $G_{n,2}$ et vérifie des conditions supplémentaires. On dira
qu'ils sont admissibles (\cite{kont} paragraphe $6.1$ et paragraphe
$2$) si
\begin{enumerate}
\item Les arêtes partent toutes des sommets de première espèce.
\item Il part deux arêtes de chaque sommet de première espèce.
\item  Le but d'une ar\^ete est diff\'erent de sa source (il n'y a pas
de boucle).
\item Il n'y a pas d'ar\^etes multiples (même source, même but).
\end{enumerate}

\begin{figure}[!h]
\begin{center}
\includegraphics[width=8cm]{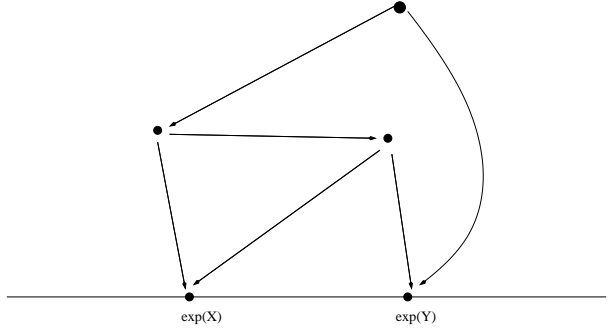}
\caption{\footnotesize Graphe admissible (de type Lie)}\label{graphe
admissible}
\end{center}
\end{figure}

\subsubsection{Opérateur différentiel associé à un graphe}

Soit $\G$ un graphe admissible de $G_{n,2}$. De tout point aérien
$i$ partent deux arêtes numérotées $(e_{i}^{a},e_{i}^{b})$.

Soit $\alpha$ un deux-vecteur sur $\R^{d}$. On peut alors associer à
tout graphe admissible un opérateur bi-différentiel sur $\R^{d}$
comme expliqué dans \cite{kont} paragraphe $2$. On notera
$B_{\G}(f,g)$ l'opérateur bidifférentiel associé que l'on suppose
agir sur les fonctions $f$ et $ g$.

Expliquons heuristiquement la formule. Sur chaque sommet aérien on
met le deux-vecteur et sur les sommets terrestres on met les
fonctions $f$ et $g$. Chaque arête arrivant sur un sommet dérive la
fonction associée au sommet. On multiplie les fonctions ainsi
obtenues et on somme sur toutes les possibilités. Concrètement la
formule est la suivante. Pour chaque arête $e$, on note par $s(e)$
le point aérien source (départ) et par $b(e)$ le point but
(arrivée). Dans la formule ci-dessous $I$ décrit l'ensemble  les
applications de l'ensemble des arêtes $E_{\G}$ dans l'ensemble des
indices de coordonnées $\{1, \cdots d\}$.

$$
B_{\Gamma,\alpha}(f,g) = \sum_{I} \bigg[\prod_{k=1}^n
\big(\prod_{\substack{e \in E_{\Gamma}\\ b(e) = k}}
\partial_{I(e)}\big) \alpha^{I(e_k^a)I(e_k^b)}\bigg]\big( \prod_{\substack{e
\in E_{\Gamma}\\b(e) = \ov{1}}}\partial_{I(e)}\big) f
\big(\prod_{\substack{e\in E_{\Gamma}\\ b(e) = \ov{2}}}
\partial_{I(e)}\big) g.$$

\subsubsection{Forme d'angles}

\begin{figure}[!h]
\begin{center}
\includegraphics[]{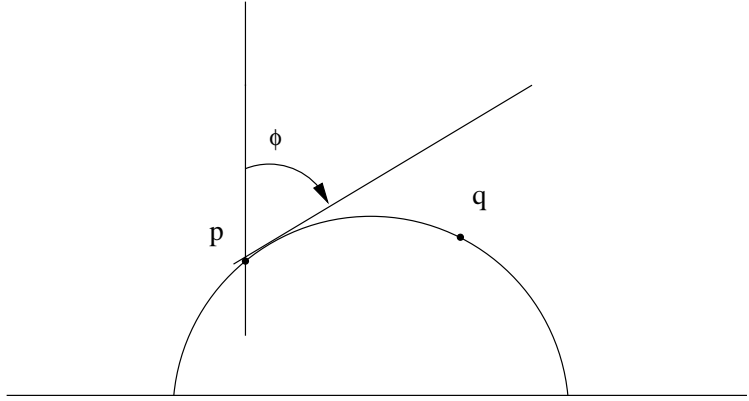}
\caption{\footnotesize Fonction d'angle}\label{angle}
\end{center}
\end{figure}

Soient deux points distincts $(p,q)$ dans le demi-plan de Poincaré
muni de la métrique de Lobachevsky. On note
\begin{equation}
\phi(p,q)=\frac{1}{2i}\log(\frac{(q-p)(\overline{q}-p)}
{(q-\overline{p})(\overline{q}-\overline{p})}).
\end{equation}
C'est l'angle entre la géodésique $(p,\infty)$ et $(p,q)$ où
l'infini peut être vu comme l'infini sur la droite réelle (figure
[\ref{angle}])

\begin{figure}[!h]
\begin{center}
\includegraphics[width=6cm]{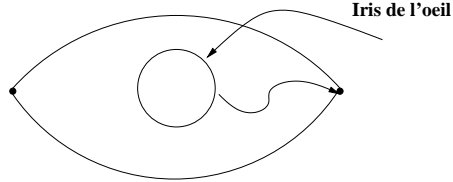}
\caption{\footnotesize La vari\'et\'e $\ov{C}_{2,0}$, dite \oe il de
Kontsevich }\label{oeil}
\end{center}
\end{figure}

La fonction d'angle s'étend  à la compactification
 $\ov{C}_{2,0}$ en une fonction d'angle régulière. La variété $\ov{C}_{2,0}$ est
 précisément décrite dans l'article de Kontsevich
 (\cite{kont} paragraphe $5.2$), c'est le fameux \oe il (figure [\ref{oeil}]).
  On remarquera, mais c'est
tautologique vu la construction des compactifications, que lorsque
les points $p, q$ s'approchent selon un angle $\theta$, la fonction
d'angle vaut précisément cet angle. Lorsque $p$ s'approche de l'axe
réel la fonction d'angle est nulle et lorsque  $q$ s'approche de
l'axe réel on obtient deux fois l'angle de demi-droite avec l'axe
réel.

Comme la fonction d'angle est régulière sur la compactification, on
peut considérer sa différentielle qui est alors une $1$-forme sur
$\overline{C}_{2,0}$.

\subsubsection{Poids associé à un graphe}

Si $\Gamma$ est un graphe admissible dans $G_{n,2}$, alors  toute
arête $e$ définit par restriction une fonction d'angle notée
$\phi_{e}$ sur la variété $\overline{C}^{+}_{n,2}$. Le produit
ordonné
\begin{equation}
\Omega_{\Gamma}=\bigwedge _{e \in E_{\Gamma}} d\phi_{e}
\end{equation}
est donc une $2n$-forme sur  $\overline{C}^{+}_{n,2}$ qui est de
dimension aussi $2n$. Le poids associé sera donc

\begin{equation}
w_{\Gamma}=\frac{1}{(2\pi)^{2n}}\int_{\overline{C}^{+}_{n,2}}
\Omega_{\Gamma}.
\end{equation}

\subsubsection{Permutation des ar\^etes}

Soit $\Gamma$ un graphe admissible dans $G_{n,2}$ et $\alpha$ un
deux-vecteur régulier sur $\R^d$. Le groupe
$\underbrace{S_{2}\times\cdots\times S_{2}}_n$, produit des groupes
de permutations des ar\^etes attach\'es \`a chaque sommet, agit
naturellement sur $\Gamma$ par permutation de l'\'etiquetage des
ar\^etes. On vérifie que l'on a ($\epsilon(\sigma)$ désigne la
signature de la permutation $\sigma$):
\begin{eqnarray*}B_{\sigma.\Gamma}&=\varepsilon(\sigma)B_\Gamma
\end{eqnarray*}
\begin{eqnarray}\label{permutation} w_{\sigma.\Gamma}&=\varepsilon(\sigma)w_\Gamma,
\end{eqnarray}
de sorte que le produit $w_\Gamma. B_\Gamma$ ne d\'epend pas de
l'\'etiquetage. On remarque aussi que le c\oe fficient $w_{\G}$ et
l'opérateur $B_{\G}$ ne dépendent pas de l'ordre des sommets.

\subsubsection{Principe de la démonstration}
L'associativité de l'étoile-produit défini par Kontsevich résulte de
deux ingré\-dients. Le premier ingrédient est la formule de Stokes
qui va donner des relations entre les c\oe fficients $w_{\G}$. Le
deuxième ingrédient est un lemme fondamental (\cite{kont} lemme
$6.6$) qui assure que les concentrations aériennes de plus de trois
points auront des contributions nulles. Au final le défaut
d'associativité de cet étoile-produit sera nul car le deux-vecteur
$\alpha$ est un deux-vecteur de Poisson (\cite{kont} paragraphe
$6.4$).

\subsection{Le cas linéaire}\label{cas lineaire}

Soit $\g$ une algèbre de Lie de dimension finie sur $\R$. L'espace
vectoriel dual $\g^*$ possède une structure de Poisson linéaire
donnée par la moitié du crochet de Lie.  Ces deux notions sont
équivalentes. Si $e_1, \cdots, e_d$ d\'esigne une base de $\g$,
$e_1^*, \cdots, e_d^*$ la base duale et $\partial_i$ la d\'eriv\'ee
dans la direction du vecteur  $e_i^*$,  alors le deux-vecteur de
Poisson associé est
$\alpha=\frac{1}{2}\sum_{i,j}[e_i,e_j]\partial_i\otimes\partial_j$.

Dans ce cas (\cite{kont} th\'eor\`eme $8.3.1$) l'ordre des
opérateurs $B_{\G}$ est suffisamment croissant de sorte que, lorsque
$f$ et $g$ sont deux fonctions polynomiales, la formule de
Kontsevich est en fait une somme finie. Cet étoile-produit vérifie
pour $X, Y$ dans $\g$ (considérés comme des fonctions linéaires sur
$\g^*$) la relation
$$X\star Y-Y\star X=h[X,Y].$$
En localisant en $h=1$, l'étoile-produit de Kontsevich définit donc
sur l'algèbre symétrique $S[\g]$ une structure isomorphe \`a
l'algèbre enveloppante de $\g$ (\cite{kont} th\'eor\`eme $8.3.1$).

\subsubsection{Géométrie des graphes}

Dans le cas linéaire, les graphes qui  vont intervenir de manière
non triviale dans la formule (\ref{formuleK}) ont une description
relativement simple. En effet chaque sommet de première espèce ne
pourra recevoir qu'au plus une arête: ce sont par d\'efinition les
graphes pertinents (relevant graphs)  et on renvoie le lecteur à
l'article \cite{AST} paragraphe $3.1$ pour une description précise
des graphes pertinents qui interviennent dans la formule finale. \\

En fait des arguments \'el\'ementaires (d\'etaill\'es dans
\cite{AST} paragraphe $3.1$) montrent que tout graphe admissible
pertinent se décompose en produit de graphes simples. Il y a deux
types de graphes simples: ceux qui contiennent une seule roue (et
donc pas de racine) comme dans la figure [\ref{roue}] on les
appellera \textit{graphes de type roue } (\cite{AST} définition
$3.1.1$) et ceux qui contiennent une seule racine (et donc pas de
roue) comme dans la figure [\ref{graphe admissible}], ce sont les
\textit{graphes de type Lie} (\cite{Ka} définition $3.1$).

\subsubsection{Symbole de $B_{\G}$}

Dans  \cite{AST} paragraphe $3$ on a associé  à chaque graphe
 admissible $\Gamma$ un symbole $a_{\G}$. C'est une fonction de
$\g\times \g$ à valeurs dans $S[\g]$   donnée par la formule pour
$X,Y$ dans $\g$

$$a_{\G}(X,Y)=B_{\G}(e^X, e^Y)e^{-X-Y}.$$

\begin{figure}[h!]
\begin{center}
\includegraphics[width=8cm]{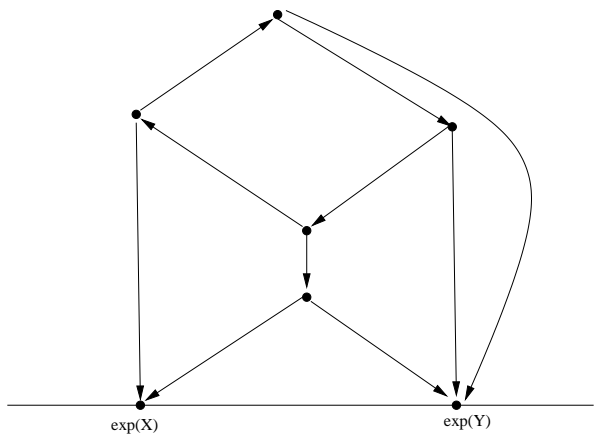}
\caption{\footnotesize Graphe simple de type Roue}\label{roue}
\end{center}
\end{figure}

Par exemple le symbole associé au graphe de la figure [\ref{graphe
admissible}] est $$a_{\G}(X,Y)=\frac{1}{8}[[X,[X,Y]],Y]$$ et le
symbole associé au graphe de la figure [\ref{roue}] est
$$a_{\G}(x,y)=\frac{1}{2^5}\tr_{\g}( \ad[X,Y]\, \ad X\, \ad Y \,\ad Y).$$

Lorsque le graphe est simple de type Lie, alors $a_{\G}(X,Y)$ est
naturellement un élément de l'algèbre de Lie engendrée par $X$ et
$Y$. Lorsque $\G$ est produit de graphes simples, le symbole est le
produit des symboles associés (\cite{AST} lemme $3.6$).

\subsubsection{Formule de Campbell-Hausdorff en termes de graphes}

Comme démontré dans \cite{Ka} th\'er\`eme $5.1$, la formule de
Campbell-Hausdorff s'écrit alors pour $X$ et $Y$ dans $\g$ (la série
est convergente pour $X$ et $Y$ proches de $0$)

\begin{equation}\label{eqCBH}
Z(X,Y)=X+Y +\sum _{n=1}^{\infty}\frac{1}{n!}\sum_{\G}w_{\Gamma}
a_{\G}(X,Y).
\end{equation}

\noindent où la somme porte sur les graphes simples de type Lie. Les
graphes qui contribuent de manière non triviale dans cette formule
n'ont donc qu'une racine et ne possèdent pas de symétries. Par
conséquent les graphes étiquetés (numérotés) associés à un graphe
géométrique (graphe orienté associé pour lequel on oublie
l'étiquetage) sont au nombre de $n!2^n$. Le terme $n!$ se compense
avec le terme  $\frac{1}{n!}$ et le terme $2^n$ disparaît aussi car
on a pris le deux-vecteur de Poisson associé à la moitié du crochet
de Lie. On notera $\Gamma(X,Y)$ le symbole associé au graphe
géométrique lorsqu'on prend le deux-vecteur
 de Poisson associé au crochet de Lie. On a alors la formule
plus synthétique suivante

\begin{equation}\label{eqCBHgeometrique}
Z(X,Y)=X+Y +\sum_{\Gamma }w_{\Gamma}\Gamma(X,Y)
\end{equation}
 où la somme porte sur les graphes géométriques simples de type Lie.
On remarquera toutefois que le symbole $\Gamma(X,Y)$ est mal défini
si le graphe n'est pas étiqueté. Pour résoudre ce problème il suffit
de remarquer que c'est aussi le cas pour le c\oe fficient
$w_{\Gamma}$ et  que ces deux difficultés se compensent grâce aux
équations (\ref{permutation}).

\subsubsection{Symbole formel de l'étoile-produit
dans le cas linéaire et fonction de densité}
 La propriété multiplicative du symbole et la
combinatoire sur les graphes assurent que le symbole formel associé
à l'étoile-produit de Kontsevich est, dans le cas des algèbres de
Lie, un produit de deux termes de type exponentiel correspondant aux
contributions des graphes simples de type Lie (c'est la formule de
Campbell-Hausdorff par \cite{Ka} théorème $5.1$) et des graphes
simples de type roue (c'est la fonction de densité par \cite{AST}
proposition $3.7$). Le symbole formel vu dans $S[\g][[h]]$, vérifie
alors (\cite{AST} proposition $3.8$)

\begin{equation}\label{eqstar}
\big(\exp(X)\star\exp(Y)\big)\exp(-X-Y)
=D(hX,hY)\exp\big(\frac{1}{h}Z(hX,hY)-X-Y\big)
\end{equation}
avec $D(X,Y)$ la fonction de densité (série convergente pour $X$ et
$Y$ proches de $0$) :
\begin{equation}
D(X,Y)=\exp\Big(\sum_{\Gamma  }w_{\Gamma}\Gamma(X,Y)\Big)
\end{equation}
où la somme porte sur l'ensemble des graphes géométriques simples et
de type roue. On montre (\cite{AST} proposition $3.12$) que cette
fonction vaut

\begin{equation}\label{eqdensite}
D(X,Y)=\frac{j^{1/2}(X)j^{1/2}(Y)}{j^{1/2}(Z(X,Y))}
\end{equation}
avec $j$ le déterminant jacobien de l'application exponentielle à
savoir
\begin{equation}\label{fonctionj}
j(X)=\det_\g\Big(\frac{1-\exp(-\ad X)}{\ad X}\Big).
\end{equation}

\subsection{Rappel sur les espaces symétriques}
Soit $\g$ une algèbre de Lie de dimension finie sur $\R$, on appelle
paire symétrique, tout couple $(\g, \sigma)$ avec $\sigma$ une
automorphisme de Lie involutif. Habituellement on note $\g=\k\oplus
\p $  la décomposition en espaces propres relativement à $\sigma$.

Lorsque $\g$ est l'algèbre de Lie d'un groupe de Lie $G$ (connexe et
simplement connexe) et $\sigma$ la différentielle d'un involution de
$G$ (un automorphisme de groupe vérifiant $\sigma^2=1$), on note
 $G^{\sigma}$ le sous-groupe des points fixes de $\sigma$ dans $G$ et
$K$ la  composante connexe de $G^{\sigma}$ qui contient l'élément
neutre $e$. L'espace $G/K$ est appelé espace symétrique. Alors $\k$
est l'algèbre de Lie de $K$, mais
 en général $\p$ n'est pas une sous-algèbre de Lie car on a
 les relations suivantes $[\k, \p] \subset \p$ et $[\p, \p]\subset \k$.
 On peut identifier
 $\p$ \`a l'espace tangent $T_{eK}(G/K)$. Il existe cependant
 sur $\p$ un produit
 triple donné par la formule pour $X,Y,Z$ dans $\p$
\begin{equation}
(X, Y, Z) \mapsto [[X, Y], Z].
\end{equation}
Les espaces vectoriels munis de produits triples sont pour les
espaces symétriques, ce que sont les algèbres de Lie pour les
groupes de Lie (\cite{loos} \textbf{vol I}, paragraphe II-$2$).

L'application exponentielle définit un difféomorphisme local de $\p$
sur $G/K$. On la note $\Exp$ et on notera $\exp_{\g}$ l'application
exponentielle de $\g$ dans $G$. On en déduit alors l'existence d'une
formule de Campbell-Hausdorff pour les espaces symétriques, définie
de la manière suivante. Pour $X$ et $Y$ dans $\p$ proches de $0$, il
existe une série convergente $Z_{sym}(X,Y)$ à valeurs dans $\p$
telle que l'on ait
$$\exp_{\g}(X)\Exp(Y)=\Exp\big(Z_{sym}(X,Y)\big).$$
En utilisant l'involution $\sigma$ on trouve facilement
\begin{equation}
\exp_{\g}(2Z_{sym}(X,Y))=\exp_{\g}(X)\exp_{\g}(2Y)\exp_{\g}(X).
\end{equation}
Cela nous suggère de modifier le crochet de Lie pour $\g$ et de
prendre plutôt deux fois le crochet. Notons $\g_2$ l'algèbre de Lie
ainsi obtenue et $G_2$ le groupe de Lie connexe et simplement
connexe associé. Notons $\exp_{\g_2}$  l'application exponentielle
associée. On aura alors (le membre de droite est calculé dans $G_2$)

\begin{equation}\label{eqsym}
\exp_{\g_2}\big(Z_{sym}(X,Y)\big)=
\exp_{\g_2}(X/2)\cdot\exp_{\g_2}(Y)\cdot\exp_{\g_2}(X/2).
\end{equation}

\section{Méthodes de déformation pour les espaces symétriques}

\subsection{Calcul de la formule de Campbell-Hausdorff}
Comment exprimer maintenant $Z_{sym}(X,Y)$ en termes de diagrammes?
Cela revient à se demander comment on calcule les produits itérés de
l'étoile-produit. La proposition suivante résout la question (voir
l'article
\cite{kont} paragraphe $8.3.3.2$ pour un calcul analogue). \\

\begin{figure}[h!]
\begin{center}
\includegraphics[width=10cm]{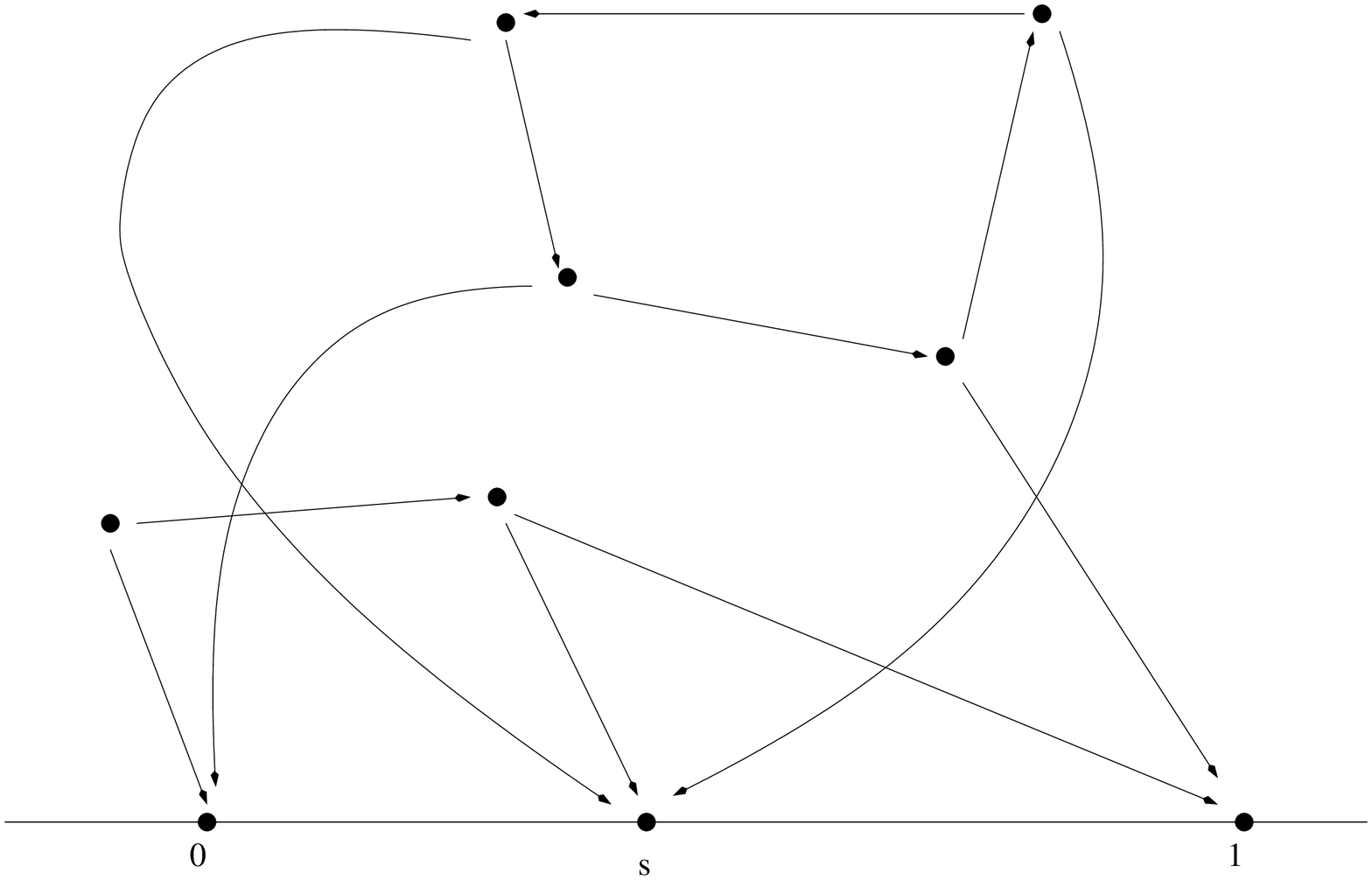}
\caption{\footnotesize Graphe avec trois points
terrestres}\label{graphes3pts}
\end{center}
\end{figure}

On considère $G_{n,3}$ l'ensemble des graphes étiquetés avec $n$
points aériens et $3$ points terrestres. On place les points
terrestres en $0, s, 1$ avec $s\in ]0,1[$ (figure
[\ref{graphes3pts}]). On étend la notion de graphes admissibles aux
graphes de $G_{n,3}$ (voir section \ref{defgrapheadmissible}).

Notons $w_{\Gamma}(s)$ la valeur du c\oe fficient obtenu par
intégration de la forme $\Omega_{\G}$ sur les configurations à $s$
fixe (en d'autres termes on intègre sur les points aériens).
On considère la structure de Poisson linéaire sur $\g^*$, donnée par
la moitié du crochet de Lie. On note $\alpha$ le deux-vecteur
associé.
\begin{prop}
Pour $f,g,h$ trois fonctions régulières sur $\g^*$ on a
$$f\star g\star k= fgk +\sum _{n=1}^{\infty}\frac{h^{n}}{n!}
\sum_{\substack{\G \in G_{n,3}\\ \G
\mathrm{admissible}}}w_{\Gamma}(s) B_{\G}(f,g,k).$$

\end{prop}
\noindent {\bf Preuve:} On doit montrer deux choses, d'une part que
l'expression de droite est indépen\-dante de $s$ et qu'elle vaut
bien le terme de gauche.

On utilise le calcul de la dérivée de $w_{\Gamma}(s)$  comme dans
\cite{To} paragraphe $4.1$. Cela revient à utiliser la formule de
Stokes  de manière infinitésimale qui est à la base de
l'associativité (\cite{kont} paragraphe $6.4$). En pratique cela
veut dire que l'on doit concentrer des points. Comme les points $0,
s, 1$ sont sur l'axe réel, les concentrations que l'on doit
effectuer sont aériennes (sinon pour des raisons de dimension on se
retrouve dans le cas \cite{kont} paragraphe $6.4.2.2$). Les
concentrations de plus de trois points aériens sont nulles
(\cite{kont} lemme $6.6$). Les seules expressions non nulles,
proviennent alors des concentrations de deux points aériens. Elles
se compensent par l'identité de Jacobi, à savoir $[\alpha,
\alpha]_S=0$ (crochet de Schouten). Le terme de droite est donc
constant.

Maintenant lorsque $s$ tend vers $0$, on obtient toutes les
concentrations de $p$ points aériens avec leurs $2p$ arêtes. En
effet, il faut tenir compte des composantes de bord dans les
compactifications des espaces de configurations (\cite{kont}
$6.4.2.1$).  Le résultat vaut alors la factorisation
$$(f\star g)\star k.$$ On en retrouve l'associativité de l'étoile-produit en
considérant la limite quand $s$ tend vers $1$.$\Box$
\\

\noindent \textbf{Remarque : } Dans la proposition ci-dessus on n'a
pas utilisé le fait que la structure de Poisson était linéaire, on
en déduit que cette proposition est vraie pour tout deux-vecteur de
Poisson.
\\

\noindent Plaçons maintenant le point médian en $s=\frac{1}{2}$. On
place en $0$ la fonction $\exp(X/2)$, en $\frac{1}{2}$ la fonction
$\exp(Y)$ et en $1$ la fonction $\exp(X/2)$. On s'intéresse à la
contribution des termes dans

$$\exp(X/2)\star\exp(Y)\star \exp(X/2)=$$
\begin{equation}\label{starSym} \exp(X+Y)+\sum
_{n=1}^{\infty}\frac{h^n}{n!}\sum_{\substack{\G \in G_{n,3}\\
\G\mathrm{admissible}} }w_{\Gamma}(\frac{1}{2})
B_{\G}\Big(\exp(X/2), \exp(Y),\exp(X/2)\Big).
\end{equation}
Le bivecteur de Poisson intervenant est comme d'habitude la moitié
du crochet de Lie (n'oublions pas que l'on travaille dans $\g_2$).
On fait maintenant intervenir la symétrie par rapport à la droite
$x=\frac{1}{2}$. Si $\G$ est un graphe, on note $\G^{\wedge}$ le
graphe symétrisé. Clairement, les quantités
$$B_{\G}\Big(\exp(X/2), \exp(Y),\exp(X/2)\Big)$$ et
$$B_{\G^{\wedge}}\Big(\exp(X/2), \exp(Y),\exp(X/2)\Big)$$ sont égales. Par
contre on a $w_{\G^{\wedge}}=(-1)^n w_{\G}$. On en déduit que les
contributions des termes avec $n$ impair sont nulles.
\\

La géométrie et la combinatoire des graphes intervenant dans la
formule (\ref{starSym}) est clairement identique à celle du
paragraphe \ref{cas lineaire}. On étend donc toutes les notions
rappelées au paragraphe \ref{cas lineaire}.
\\

On note de la même façon que précédemment le symbole de l'opérateur
$B_{\G}$ lorsqu'on prend pour deux-vecteur de Poisson le crochet de
Lie de $\g_2$. On a donc
\begin{equation}\label{symbolesym}
\G(X,Y)=B_{\G}\Big(\exp(X/2), \exp(Y),\exp(X/2)\Big)\exp(-X-Y).
\end{equation}
On en déduit, comme dans l'équation (\ref{eqstar}), que le membre de
gauche dans (\ref{starSym}) est un produit de deux termes de type
exponentiel correspondant aux contributions des graphes simples de
type roue et aux contributions des graphes simples de type Lie. La
contribution dans (\ref{starSym})  des graphes admissibles dont les
composantes simples sont de type Lie calcule le terme de droite dans
(\ref{eqsym}) (c'est le théorème $5.1$ de \cite{Ka}). Par conséquent
la contribution des graphes simples de type Lie dans (\ref{starSym})
correspond à la formule de Campbell-Hausdorff pour les espaces
symétriques par (\ref{eqsym}).
\\

En considérant les graphes géométriques associés, on en déduit comme
pour la formule (\ref{eqCBHgeometrique}) la proposition suivante :
\begin{prop}
La formule de Campbell-Hausdorff pour les espaces symétriques
s'écrit pour $X$ et $Y$ dans $\p$ :
$$Z_{sym}(X,Y)=X+Y+ \sum_{\G} w_{\Gamma}(\frac{1}{2}) \G(X,Y)$$
avec la somme portant sur les graphes simples géométriques de type
Lie avec trois points terrestres. En particulier on a bien
$Z_{sym}(X,Y) \in \p$.
\end{prop}

Pour les mêmes raisons que \cite{ADS} lemme $2.2$ et $2.3$ ou
\cite{AST} proposition $3.9$, les séries associées sont convergentes
pour $X$ et $Y$ dans $\p$ proches de $0$. On peut donc prendre $h=1$
et on obtient
\begin{equation}
\exp(X/2)\star\exp(Y)\star
\exp(X/2)=D_{sym}(X,Y)\exp\big(Z_{sym}(X,Y)\big).
\end{equation}
Notons $J_{\g}(X)$ le jacobien de l'application $\Exp$. On a la
formule classique
\begin{equation}
J_{\g}(X)=\det_\p\Big(\frac{\sinh \ad X}{\ad X}\Big).
\end{equation}
En utilisant la fonction de densité (\ref{eqdensite}) pour les
algèbres de Lie  et l'associativité de l'étoile-produit, on peut
calculer sans difficulté la fonction de densité $D_{sym}(X,Y)$. On
obtient
\begin{equation}\label{eqdensitesym}
D_{sym}(X,Y)=\frac{j_{\g_2}(X/2)j_{\g_2}(Y)^{1/2}}{j_{\g_2}(Z_{sym}(X,Y))^{1/2}}
\end{equation}
où la fonction $j_{\g_2}$ (\ref{fonctionj}) est calculée pour le
crochet double. Or pour $X\in \p$ on a
$$j_{\g_2}(X)=\exp(-\tr_{\g}(\ad X)) J_{\g}(X)^2$$ avec membre de
droite calculé pour le crochet simple. On en déduit le lemme
\begin{lem}
On a la formule
\begin{equation}\label{eqdensitesymbis}
D_{sym}(X,Y)=\frac{J_{\g}(X/2)^2 J_{\g}(Y)}{J_{\g}(Z_{sym}(X,Y))}.
\end{equation}
\end{lem}
\subsection{Première Déformation}
On va effectuer maintenant une première déformation de la formule de
Camp\-bell-Hausdorff et de la fonction de densité. Elle consiste à
déformer les c\oe fficients $w_{\G}$ en déplaçant les points
terrestres dans le demi-plan de Poincaré (\cite{kont} $8.2.3$).
\\

\begin{figure}[h!]
\begin{center}
\includegraphics[width=10cm]{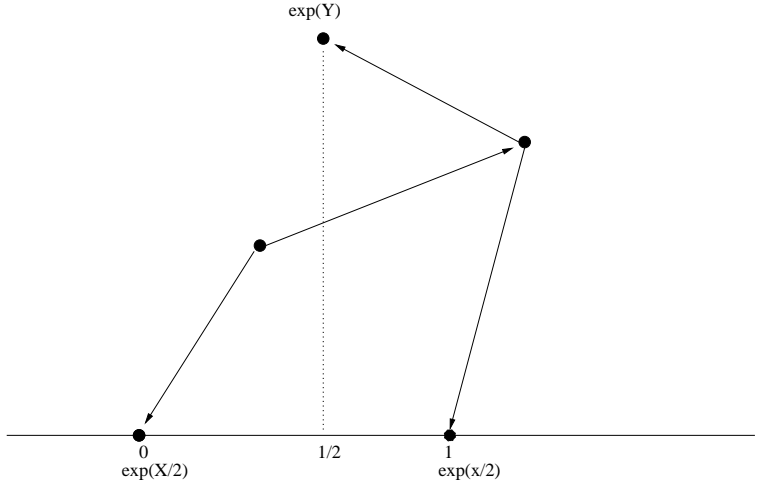}
\caption{\footnotesize Première déformation typique
}\label{premiereDeformation}
\end{center}
\end{figure}

La première déformation correspond typiquement au graphe de la
figure [\ref{premiereDeformation}] où la position du point associé à
$\exp(Y)$ est repéré par la variable $t$. Cette variable désigne la
position du point dans $\ov{C_{1,2}}^+$, avec le point aérien sur la
droite d'équation $x=\frac{1}{2}$. Cela revient à considérer des
graphes simples de type Lie de $G_{n,3}$, où on a placé les points
de deuxième espèce en $(0,t,1)$.

On ajoute un indice $t$ pour signifier la dépendance en $t$. La
règle concernant les symétries dans les c\oe fficients $w_{\G}(t)$
joue de la même façon que précédemment. Seuls les termes avec un
nombre pair de sommets de première espèce vont intervenir.

\begin{figure}[h!]
\begin{center}
\includegraphics[width=10cm]{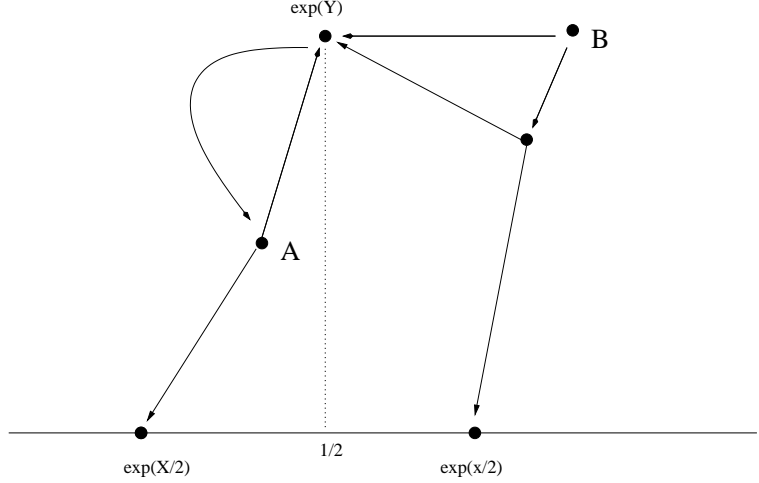}
\caption{\footnotesize Contribution Type dans la
dérivée}\label{coefderive}
\end{center}
\end{figure}

On effectue comme dans \cite{To} paragraphe $4.1$ le calcul de la
dérivée de la fonction
\begin{equation} Z_{sym,t}(X,Y)= X+Y+ \sum_{\G}
w_{\Gamma}(t) \G(X,Y). \end{equation} On ne doit considérer que les
concentrations d'un point aérien sur $\exp(Y)$, les autres
concentrations donnant soit directement $0$ (\cite{kont} lemme
$6.6$) soit s'annulant par compensation due à l'identité de Jacobi.
On se retrouve avec des c\oe fficients dérivés correspondant à des
graphes du genre de la figure [\ref{coefderive}].

Le c\oe fficient $w_{A,B}(t)$ correspondant à un graphe comme dans
la figure [\ref{coefderive}] contient un terme en $dt$. Pour les
mêmes raisons que celles expliquées dans \cite{To} proposition
$4.1$, le c\oe fficient $w_{A,B}(t)$ se factorise en un produit
$\widetilde{w_{A}}(t)w_{B}(t)$ où $w_B(t)$ se calcule à $t$ constant
et $\widetilde{w_A}(t)$ contient toute la contribution de la
différentielle en $t$. En regroupant les graphes qui dégénèrent de
la même manière, on fait apparaître un terme, au niveau des
symboles,  du genre
\begin{equation}
\partial_{Y}B(X, Y)\cdot [Y, A(X,Y)]
\end{equation}
où $A, B$ désignent par abus les symboles des graphes $A, B$. Ces
graphes sont forcément des graphes simples admissibles de type Lie
(avec points de deuxième espèce placés en $(0,t,1)$). Remarquons par
ailleurs qu'il s'agit bien du crochet de Lie, car on a pris le
deux-vecteur de Poisson associé au crochet de Lie. Pour des raisons
de symétrie, on aura une contribution non triviale des c\oe
fficients lorsque le graphe $B$ a un nombre pair de sommets de
première espèce et le graphe $A$ un nombre impair de sommets  de
première espèce (à cause de l'arête qui joint le point maintenant
aérien $t$ à la racine de $A$). La factorisation du c\oe fficient
permet alors de tout factoriser comme dans \cite{To} théorème $4.2$
:

\begin{equation}\sum_{A,B} w_{A,B}(t) [Y,A(X,Y)]\cdot
\partial_{Y}B(X,Y)=\sum_{A,B}w_{B}(t)\widetilde{w_{A}}(t)[Y,A(X,Y)]\cdot
\partial_{Y}B(X,Y)
\end{equation}
\begin{equation}
[Y, \sum_{A}\widetilde{w_{A}}(t)A(X,Y)]\cdot
\partial_{Y} \left(\sum_{B} w_{B}(t) B(X,Y)\right)=[Y,F_t]\cdot
\partial_{Y}Z_{sym,t}(X,Y). \end{equation}

La convention des signes de \cite{AMM} $I.2.1$, pour l'orientation
des strates donne le signe final. Par ailleurs l'étiquetage des
graphes n'a plus d'importance une fois qu'on a considéré le produit
du c\oe fficient et du symbole (car il y a compensation des signes).
On peut donc travailler directement avec les graphes géométriques.

Pour la fonction de densité on procède comme dans \cite{To}
paragraphe $5.2$ et on obtient alors les mêmes équations que celles
que nous avons établies dans \cite{To} théorème $4.2$ et théorème
$5.2$. Par suite on peut énoncer le théorème suivant :

\begin{theo}\label{equationdiff}Soit $(\g, \sigma)$ une paire
symétrique réelle. Il existe une série convergente en $X, Y$ dans un
voisinage de $0$ dans  $\p$, on la note $F_{t}(X,Y)$. Elle est à
valeurs dans $\k$ et c'est une $1$-forme en $t$. Il existe une
déformation de la formule de Campbell-Hausdorff pour les espaces
symétriques, on la note $Z_{sym, t}(X,Y)$ et  il existe enfin une
déformation de la fonction de densité, on la note $D_{sym, t}(X,Y)$
telles  que l'on ait
\begin{equation}\label{eqdiffCBH}{ \mathrm d}Z_{sym,
t}(X,Y)=[Y,F_{t}(X,Y)]\cdot
\partial_{Y}Z_{sym, t}(X,Y)
\end{equation}
 et
\begin{equation}\label{eqdiffdensite}{ \mathrm d}D_{sym, t}(X,Y)=[Y,F_{t}(X,Y)]\cdot \partial_{Y}D_{sym,
t}(X,Y)+ \tr_{\g}\Big(D_{Y}F_t \circ \ad (Y)\Big)D_{sym,
t}(X,Y).\end{equation}
\end{theo}
\noindent La $1$-forme $F_{t}(X,Y)$, correspond à
\begin{equation} \sum_{A}\widetilde{w_A}
(t)A(X,Y) \end{equation}
 avec $A$ un graphe simple géométrique de
type Lie (où les points de deuxième espèce sont placés en $(0,
t,1)$) et $w_A$ le c\oe fficient calculé en ajoutant une arête du
point $t$ vers le sommet de $A$ comme dans la figure
[\ref{coefderive}]. La notation $D_Y F_t(X,Y)$ signifie que l'on a
pris la différentielle partielle de $F_t$ en $Y$.

\subsection{Modification de la fonction de densité pour les espaces
\\ symétriques: cas résoluble et très symétrique} Les équations
obtenues dans le théorème précédent ne sont pas satisfaisantes du
point de vue des espaces symétriques. En effet au lieu de trouver
une $\tr_{\g}$ on aurait aimé trouver une $\tr_{\k}$ dans
(\ref{eqdiffdensite}).

Comme $D_{Y}F_t$ et $\ad (Y)$ sont des endomorphismes de $\g$ qui
envoient $\k$ sur $\p$ et $\p$ sur $\k$, on en déduit que l'on a
\begin{equation}
\tr_{\p}\Big(\ad(Y)\circ D_{Y}F_t\Big)=\tr_{\k}\Big(D_{Y}F_t \circ
\ad (Y)\Big)
\end{equation}
puis \begin{equation}\label{tracesurp}\tr_{\g}\Big(D_{Y}F_t \circ
\ad (Y)\Big)=2\tr_{\p}\Big(\ad(Y)\circ
D_{Y}F_t\Big)-\tr_{\k}\Big([D_{Y}F_t, \ad (Y)]\Big)
\end{equation}
\begin{lem}[\cite{rou86} page $573$]\label{lemmerou}
Pour les espaces symétriques résolubles et les paires très
symétriques le commutateur $[D_{Y}F_t, \ad (Y)]$  est \`a trace
nulle sur $\k$.
\end{lem}
On en déduit alors, que l'on peut transformer facilement les traces
sur $\g$ en des traces sur $\p$ à condition que le terme
$\tr_{\k}([D_{Y}F_t, \ad Y])$ soit nul. C'est exactement l'argument
de Rouvière. En particulier on trouve que ce terme est nul dans le
cas des paires symétriques résolubles et des paires très symétriques
(\cite{rou86} paragraphe $5.3$). Une paire très symétrique est la
donnée d'une paire symétrique et d'un automorphisme commutant aux
$\ad X $ et qui envoie $\p$ sur $\k $ et $\k$ sur $\p$ (\cite{rou86}
paragraphe $5$). Pour les paires très symétriques on a facilement
(\cite{rou86} paragraphe $5.3$ page $573$)
\begin{equation}
\tr_{\k}\Big(D_{Y}F_t \circ \ad Y\Big)=\tr_{\p}\Big(D_{Y}F_t \circ
\ad Y\Big)=\tr_{\k}\Big(\ad Y \circ D_{Y}F_t\Big).
\end{equation}
Des exemples connus de paires très symétriques, sont les paires
$\g\times\g/\g$ et $\g_{\C}/\g_{\R}$. Dans le premier cas, on prend
$(X,Y) \mapsto (X, -Y)$ et dans le second cas on prend la
multiplication par $i$.
\\

\noindent {\bf Conclusion: } Dans le cas des espaces symétriques il
est utile de prendre la racine carrée de la fonction de densité pour
obtenir une trace sur $\k$ (afin de compenser le c\oe fficient $2$
dans (\ref{tracesurp})).
\subsection{Calcul de la fonction de densité à l'infini}
On fait tendre $t$ vers l'infini, ce qui revient au même de fixer le
point aérien en $\frac{1}{2} +i$ et de faire tendre les deux points
terrestres vers $\frac{1}{2}$. On détermine facilement la valeur de
la fonction de densité à l'infini.
\\

\begin{figure}[h!]
\begin{center}
\includegraphics[width=10cm]{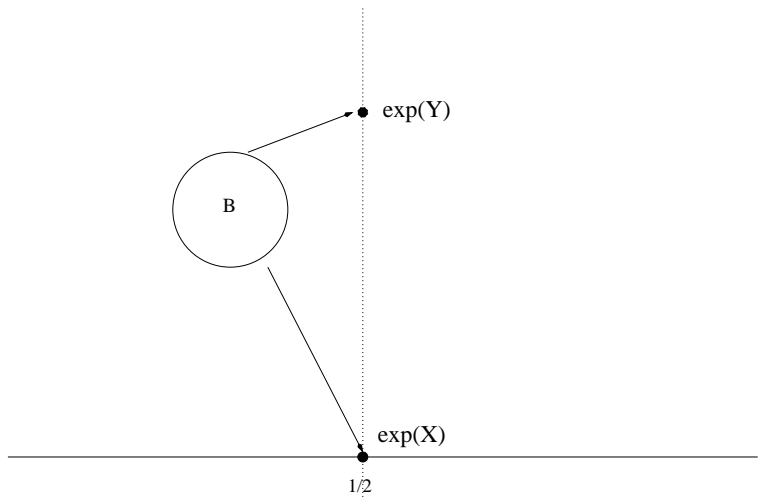}
\caption{\footnotesize Contribution Type dans la fonction de densité
$E_{1/2}(X,Y)$ }\label{densite2}
\end{center}
\end{figure}

Avant cela introduisons une autre fonction de densité que l'on note
$E_{1/2}(X,Y)$,  correspondant aux contributions des graphes
admissibles dont les composantes simples sont de type roue  et dont
les points de deuxième espèce sont placés en $\frac{1}{2}$
(correspondant à $\exp(X)$) et en $\frac{1}{2} +i$ (correspondant à
$\exp(Y)$) (voir figure [\ref{densite2}]). Par ailleurs les points
de première espèce sont associés au deux-vecteur de Poisson associé
au crochet de Lie.

\begin{figure}[h!]
\begin{center}
\includegraphics[width=10cm]{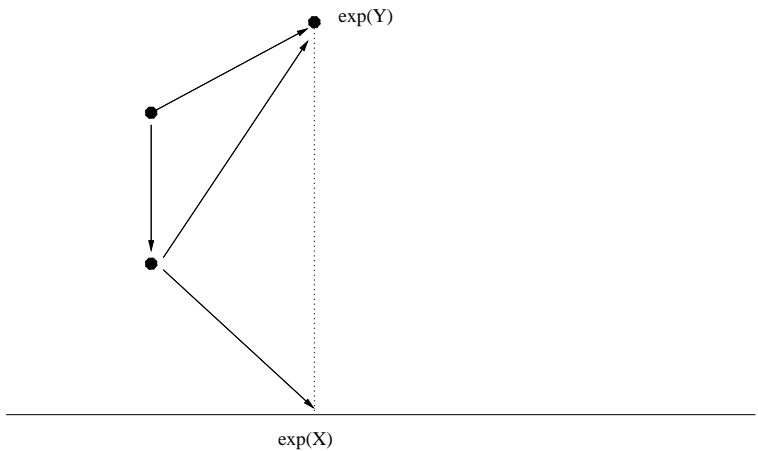}
\caption{\footnotesize Contribution Type dans $Z_{1/2}$
}\label{dessinsym2}
\end{center}
\end{figure}

De manière analogue on notera $\exp (Z_{1/2}(X,Y))$ la contribution
des graphes admissibles dont les composantes simples sont de type
Lie et dont les points de deuxième espèce sont placés en
$\frac{1}{2}$  et en $\frac{1}{2} +i$ (voir figure
[\ref{dessinsym2}]).

\begin{prop}
 Lorsque le point $t$ tend vers l'infini dans la direction imaginaire, on obtient
\begin{equation}\label{eqdensiteinfini}
D_{sym, \infty}(X,Y)=E_{1/2}(X,Y)\frac{J_{\g}^2(X/2)}{J_{\g}(X)}
\end{equation}
\end{prop}
\noindent {\bf Preuve: } En approchant les deux points terrestres
vers $\frac{1}{2}$, on doit tenir compte de toutes les
concentrations de points aériens sur ces points terrestres. Cela
ajoute un terme correspondant à la fonction de densité pour les
algèbres de Lie (\ref{eqdensite}), où l'on place $\exp(X/2)$ et
$\exp(X/2)$ aux points terrestres, ce qui donne
$$\frac{j_{\g_2}(X/2)}{j_{\g_2}(X)^{1/2}}=\frac{J_{\g}^2(X/2)}{J_{\g}(X)}. \Box $$

\begin{figure}[h!]
\begin{center}
\includegraphics[width=6cm]{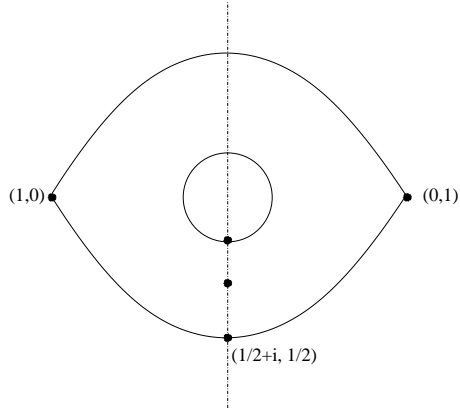}
\caption{\footnotesize Le point de déformation reste sur l'axe
$x=\frac{1}{2}$ }\label{paupiereMoitie}
\end{center}
\end{figure}

\subsection{Interprétation dans l'\oe il de Kontsevich}
Lorsqu'on a rapproché nos deux points terrestres en $\frac{1}{2}$,
on retrouve des graphes admissibles où les points de seconde espèce
sont placés en  $(\frac{1}{2} +i,\frac{1}{2})$. C'est donc un point
dans $\ov{C_{2, 0}}$ positionné sur la paupière de l'\oe il et sur
l'axe de symétrie voir figure ([\ref{paupiereMoitie}]). Cela
explique notre choix de notation : on est à la moitié sur la
paupière.

\subsection{Application à la convolution des distributions
$K$-invariantes sur les espaces symétriques résolubles ou très
symétriques}

En reprenant d'une part les arguments énoncés dans \cite{AST}
paragraphe $4$ et \cite{KV} paragraphe $2$ sur la convolution des
germes de distributions invariantes et  d'autre part les définitions
de \cite{rou86} paragraphe $3$ pour la convolution  sur $\p$, on
déduit le théorème suivant (pour simplifier les notations on a noté
$J$ pour $J_{\g}$) :

\begin{theo}\label{theoconv}
Soit $(\g, \sigma)$ une paire symétrique résoluble ou très
symétrique. Soient $u(x)$ et $v(y)$ des germes de distributions $\ad
\k$-invariantes en $0$ dans $\p$ à supports compatibles. On a alors
la formule suivante \begin{eqnarray*}\left\langle u(x)v(y),
\frac{J^{1/2}(x)J^{1/2}(y)}{J^{1/2}(z_{sym}(x,y))}
\Phi(z_{sym}(x,y))\right\rangle = \left\langle u(x)v(y),
E_{1/2}^{1/2}(x,y)\Phi(z_{1/2}(x,y))\right\rangle
\end{eqnarray*}
\noindent pour $\Phi$ une fonction test à support dans un voisinage
de $0$ dans $\p$.
\end{theo}
\noindent {\bf Preuve:} Les problèmes de convergence et la
définition de supports compatibles sont étudiés dans \cite{rou86}
paragraphe $3.3$. La méthode est par ailleurs strictement identique.
On se contente donc de la partie formelle.

On a introduit une déformation de la formule de Campbell-Hausdorff
et de la fonction de densité. On a donc une dépendance en $t$
\begin{equation}
\Psi(t)=\left\langle u(x)v(y), D_{sym,
t}^{1/2}(x,y)\Phi(z_{sym,t}(x,y)) \right\rangle  .
\end{equation}
On calcule facilement la dérivée en $t$, compte tenu des équations
(\ref{eqdiffCBH}) et (\ref{eqdiffdensite}) du théorème
\ref{equationdiff}. L'invariance du germe $v(y)$ signifie que l'on a
\begin{equation}
v(y)[y,F_{t}(x,y)]\cdot \partial_y=-\tr_{\p}\Big( \ad(y)\circ
D_{Y}F_t\Big)v(y),
\end{equation}

On obtient alors
\begin{equation}
\Psi'(t)=\left\langle u(x)v(y), m_t(x,y)D_{sym,
t}^{1/2}(x,y)\Phi(z_{sym,t}(x,y)) \right\rangle
\end{equation}
avec $m_t(x,y)= -\frac{1}{2}\tr_{\k}([D_{Y}F_t, \ad y])$. En effet
l'invariance fait disparaître le terme en $\tr_{\p}$ dans
(\ref{tracesurp}) (n'oublions pas que l'on a pris la racine carrée).
Lorsque $(\g,\sigma)$ est une paire très symétrique ou une paire
résoluble, la  trace sur $\k$ du crochet est nulle par le lemme
\ref{lemmerou}. On en déduit que la fonction $\Psi$ est constante.

Il vient $\Psi(0)=\Psi(\infty)$, puis en utilisant
(\ref{eqdensitesymbis}) et (\ref{eqdensiteinfini})
\begin{eqnarray}
\left\langle u(x)v(y), \frac{J(x/2)
J(y)^{1/2}}{J(z_{sym}(x,y))^{1/2}}\Phi(z_{sym}(x,y))\right\rangle = \nonumber\\
\left\langle u(x)v(y),
E_{1/2}^{1/2}(x,y)\frac{J(x/2)}{J(x)^{1/2}}\Phi(z_{1/2}(x,y))
\right\rangle  .
\end{eqnarray}
En changeant un peu l'ordre des termes et en simplifiant par
$J(x/2)$ (la fonction $J$ est invariante par $K$ et l'énoncé vaut
pour tous germes de  distributions invariantes) on trouve
\begin{eqnarray}
\left\langle u(x)v(y),
\frac{J(x)^{1/2}J(y)^{1/2}}{J(z_{sym}(x,y))^{1/2}}\Phi(z_{sym}(x,y))
\right\rangle= \nonumber\\ \left\langle u(x)v(y),
E_{1/2}^{1/2}(x,y)\Phi(z_{1/2}(x,y)) \right\rangle  . \hfill
\Box\end{eqnarray}

\noindent {\bf Conclusion:} dans les cas de paires symétriques
étudiées, on a trouvé une déformation adéquate vis à vis de la
convolution des germes de distributions invariantes sur les espaces
symétriques.

\subsection{Deuxième déformation}

On utilise maintenant la deuxième déformation, qui consiste à
déplacer dans $\ov{C_{2,0}}$ notre point de la paupière vers l'iris
de l'\oe il tout en restant sur l'axe de symétrie vertical de l'\oe
il. En d'autres termes, on fait se rapprocher les deux points  de
seconde espèce  sur l'axe $x=\frac{1}{2}$.
\\

C'est la déformation que l'on a utilisé dans \cite{To}. La situation
est alors stricto sensu identique. Les équations que l'on obtient
sont celles de \cite{To} théorème $4.2$ et théorème $5.2$. Compte
tenu des symétries pour les c\oe fficients et le fait que $X$ et $Y$
sont dans $\p$, il est certain que nos champs adjoints sont dans
$\k$. Comme on restreint notre étude au cas des espaces symétriques
résolubles ou très symétriques, on peut grâce au lemme
\ref{lemmerou}, transformer dans l'équation différentielle
(\ref{eqdiffdensite}) la trace sur $\g$ en une trace sur $\k$ (ou
sur $\p$) à condition de prendre la racine carrée de la fonction de
densité.

Sur l'iris de l'\oe il cette fonction de densité vaut $1$ d'après
\cite{sh} et la déformation de CBH vérifie $Z_{1/2, iris}(X,Y)=X+Y$.
On a donc établi le théorème suivant (on note $J$ pour $J_{\g}$):

\begin{theo}\label{theoconvfinal}
Soit $(\g, \sigma)$ une paire symétrique résoluble ou très
symétrique. Soient $u(x)$ et $v(y)$ des germes de distributions $\ad
\k$-invariantes en $0$ dans $\p$ à supports compatibles. Pour $\Phi$
une fonction test à support dans un voisinage de $0$ dans $\p$, on a
alors la formule suivante
$$\left\langle u(x)v(y),
\frac{J^{1/2}(x)J^{1/2}(y)}{J^{1/2}(z_{sym}(x,y))}
\Phi(z_{sym}(x,y))dx dy= \int u(x)v(y)\Phi(x+y) \right\rangle  . $$
\end{theo}

\noindent {\bf Remarques finales:} \begin{enumerate}\item Quand
$(\g, \sigma)$ est une paire symétrique résoluble, le théorème
ci-dessus est démontré par Rouvière \cite{rou86} théorème $3.6$ en
utilisant une déformation homogène de la formule de
Campbell-Hausdorff comme dans \cite{KV}.
\item Les déformations utilisées dans cet article sont régulières
mais leur jonction n'est pas différentiable. Cela explique notre
choix d'exposition en deux parties et nous permet heureusement de
retrouver les résultats de \cite{rou86}.

\item Lorsque $u(x)$ et $v(y)$ correspondent à des distributions de
support $0$, on retrouve que la symétrisation modifiée par la racine
carrée du jacobien, constitue un isomorphisme d'algèbre de
$S[\p]^{\k}$ (les symboles $K$-invariants) sur $(U(\g)/U(\g)\cdot
\k)^{\k}$, l'algèbre des opérateurs différentiels invariants sur
l'espace symétrique $G/K$. C'est l'isomorphisme de Rouvière dans le
cas des espaces symétriques résolubles (voir \cite{hab} pour une
synthèse sur ce problème).

\item On sait qu'en général on ne peut pas obtenir un tel transport des
distributions invariantes en utilisant une modification aussi simple
de l'application exponentielle. On peut s'en convaincre par exemple
en examinant ce qu'est l'homomorphisme d'Harish-Chandra du cas
semi-simple (non très symétrique). Pour pallier cette difficulté
Rouvière a introduit une fonction auxiliaire dans \cite{rou90},
\cite{rou91}, \cite{rou94} pour les fibrés en droites sur $G/K$,
notée $e_{\lambda}(x,y)$. Il semble opportun de relier ces résultats
avec les méthodes décrites dans cet article. Ce programme de
recherche est annoncé dans \cite{hab}.

\end{enumerate}

%
%

\vskip 12pt
\parindent=0pc
C. Torossian\\
UMR 8553 du CNRS, DMA, École Normale Supérieure\\
45 rue d'Ulm 75230 Paris cedex 05
\end{document}